\documentclass{conm-p-l}

\usepackage{amsmath,amsthm,amssymb}
\usepackage{diagrams}
\usepackage{graphicx}

\newtheorem{Theorem}{Theorem}[section]

\theoremstyle{definition}
\newtheorem{Definition}[Theorem]{Definition}

\newcommand{\Set}{\mathbf{Set}}

\newcommand{\Quan}{\mathbf{Quan}}
\newcommand{\KDiag}{\mathbf{KDiag}}

\newcommand{\Mat}{\mathbf{Mat}}

\newcommand{\Span}{\mathbf{Span}}
\newcommand{\CoSpan}{\mathbf{CoSpan}}

\newcommand{\Tang}{\mathcal{T}ang}

\begin{document}

\title{Categorifying Coloring Numbers}

\author{John Armstrong}
\address{Department of Mathematics, Tulane University, New Orleans, Louisiana 70118}
\email{jarmstro@tulane.edu}
\thanks{I am deeply indebted to the input and advice of John Baez and J. Scott Carter on the preliminary versions of this paper, and to Sam Lomonaco and Louis Kauffman in the development of these ideas.}

\date{March 11, 2008.}


\begin{abstract}
Coloring numbers are one of the simplest combinatorial invariants of knots and links to describe.  And with Joyce's introduction of quandles, we can understand them more algebraically.  But can we extend these invariants to tangles -- knots and links with free ends?  Indeed we can, once we categorify.

Starting from the definition of coloring numbers, we will categorify them and establish this extension to tangles.  Then, decategorifying will leave us with matrix representations of the monoidal category of tangles.
\end{abstract}

\maketitle

\section{Introduction}
\subsection{Topological Quantum Computation and Tangle Representations}
The rise of topological quantum computation as a method to provide fault-tolerance for quantum computers\cite{MR1951039,MR1943131,MR1750541} brings with it the need to turn knot theory into representation theory.  Every computation is actually approximating a topological invariant of the knotted paths anyons follow, and every knot invariant should give a quantum computer.

But we cannot simply consider these as invariants of knots.  Computations take place through time, and we must be able to understand what happens in the first half of a computation as separate from what happens in the second half.  When we consider less than the complete run of a topological quantum computer we do not find neatly knotted paths of anyons, but rather a loose collection of tangled paths with free ends hanging out at the beginning and end of the computation.

Thus we must consider tangles\cite{MR1020583} as a natural generalization of knots and links, and a simpler one for the purposes of topological quantum computation.  To describe a topological quantum computer corresponding to a tangle we must select one transition matrix for each of four simple generating tangles, subject to a short list of conditions.  That is, we must define a matrix representation of the category $\Tang$ of tangles.  This specifies not only the evolution of the computer's state as we move anyons around each other, but also the initial conditions and the measurements to be performed as we pair them off.

And so quantum computation requires us to consider the representation theory of tangles, and to think of knot invariants as restrictions of these representations.

\subsection{Colorings and Quandles}
In this paper we will lay out this picture for a particularly simple combinatorial invariant of knots and links: the number of colorings of a link by a given involutory quandle.

Colorings of knot and link diagrams go back to Fox\cite{MR0445489}, who asked if we can color the arcs of a diagram red, green, and blue, so that at each crossing either one color appears or all three do.  More generally, in how many ways can we manage this?  It turns out that this number of colorings depends only on the knot type and not on the particular diagram.

Quandles were introduced by Joyce\cite{MR638121} and Matveev\cite{MR672410} (under the name, ``distributive groupoids'') as a tool for studying oriented knots and links.  The special case of involutory quandles first appeared as ``keis''\cite{Takasai-Keis}.  These fill the same role for unoriented links that general quandles do for oriented links.

The connection between quandles and colorings is that a coloring is essentially a homomorphism of quandles\cite{MR1675756}.  First, the set of colors $\{\mathrm{red, green, blue}\}$ can be given the structure of an involutory quandle.  Then, to an unoriented link diagram we can assign a ``fundamental'' involutory quandle encoding exactly those relations demanded by the diagram's crossings.  Fox's colorings, then, are homomorphisms from this fundamental involutory quandle to the quandle of colors.  Replacing this target quandle with other involutory quandles gives a rich stock of invariants to investigate.

The framework of quandles has been extended to include a cohomology theory analogous to that of groups\cite{MR1725613, carter-2006}.  Link colorings by various sorts of quandles have also been extensively studied\cite{MR1483791, MR1634462, MR1634467, MR1905694, MR2240912, MR2275096, MR2346499}.  However, these invariants must be extended to tangles for our purposes!

Our first step will be to ``categorify'' the coloring number invariants by considering instead the set of colorings of a given diagram.  It is essential at this point to note that this set is \emph{not} invariant under the Reidemeister moves -- only its cardinality is.  This leads us in passing to define it as an example of a link (or tangle) ``covariant''.

Next we extend our definition to cover tangles by introducing the category of spans, as defined by B\'enabou\cite{MR0220789}.  We find that defining the colorings of a tangle to be a span of sets gives us exactly the handles we need to compose them properly, and to define colorings as a functor on the category of tangles.

Finally, we ``decategorify'' our spans to find matrices\cite{Baez-QG-spring2004}.  This gives us our sought-after matrix representation of the category of tangles.  When we regard a link as a tangle, our representation will give us a $1\times1$ matrix whose single entry is the old number of colorings.

\subsection*{Acknowledgements}
I am deeply indebted to the input and advice of John Baez and J. Scott Carter on the preliminary versions of this paper, and to Sam Lomonaco and Louis Kauffman in the development of these ideas.

\section{Quandle Coloring Numbers}
\subsection{Quandles}
A ``quandle'' is an algebraic structure consisting of a set $Q$ and two binary operations $\triangleright$ and $\triangleleft$.  These satisfy the three conditions
\begin{list}{}{}
\item[\textbf{Q1}.] For all $a\in Q$, $a\triangleright a=a$.
\item[\textbf{Q2}.] For all $a,b\in Q$, $(b\triangleright a)\triangleleft b = a = b\triangleright(a\triangleleft b)$.
\item[\textbf{Q3}.] For all $a,b,c\in Q$, $a\triangleright(b\triangleright c)=(a\triangleright b)\triangleright(a\triangleright c)$.
\end{list}

As is usual for algebraic structures, we have a notion of a ``quandle homomorphism'' $f:Q_1\to Q_2$, which is simply a function from the underlying set of $Q_1$ to that of $Q_2$ which preserves the two quandle operations.  We then have the category $\Quan$ of quandles and quandle homomorphisms, which will feature prominently in our discussion.

It is useful to keep the following quandles in mind as examples.

Given any group $G$, the conjugation  $\mathrm{Conj}(G)$ with the same underlying set as $G$.  We define the operations by conjugation within the group:
\begin{align*}
b\triangleright a &= bab^{-1}\\
a\triangleleft b &= b^{-1}ab\\
\end{align*}

If $G$ is abelian, then the operations in $\mathrm{Conj}(G)$ are trivial.  But we do have another interesting quandle structure.  The dihedral quandle $D(G)$ also has the same underlying set as $G$, but we now define the two operations:
\begin{equation*}
b\triangleright a = 2b-a = a\triangleleft b
\end{equation*}
This quandle satisfies an additional condition
\begin{list}{}{}
\item[\textbf{QInv}.] For all $a,b\in Q$, $b\triangleright a = a\triangleleft b$
\end{list}
When this condition is satisfied, we say the quandle is ``involutory''.

\subsection{Colorings}
\begin{figure}[hbt]
  \centerline{\includegraphics{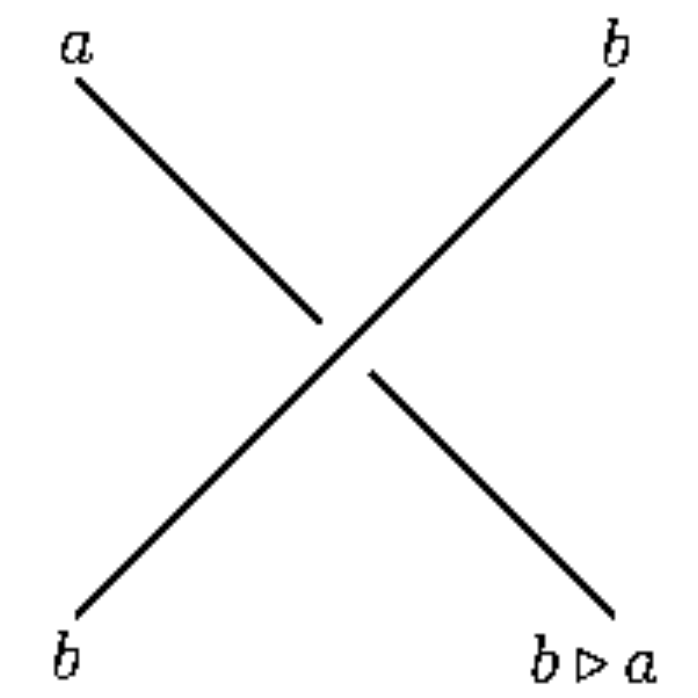}}
\caption{Coloring arcs at a crossing}\label{fig:CrossingColor}
\end{figure}

Given a unoriented knot or link diagram and an involutory quandle $X$, we color the diagram by assigning an element of $X$ to each arc of the diagram.  When an arc with color $a$ meets an overcrossing arc with color $b$, the arc on the other side must be colored $b\triangleright a$, as in figure \ref{fig:CrossingColor}.

Notice here that it doesn't matter which undercrossing arc we regard as coming in and which we regard as going out of the crossing because we are using an involutory quandle.  The axioms \textbf{QInv} and \textbf{Q2} tell us that
\begin{equation*}
b\triangleright(b\triangleright a)=b\triangleright(b\triangleleft a)=a
\end{equation*}

\begin{figure}[p]
  \includegraphics[width=\textwidth]{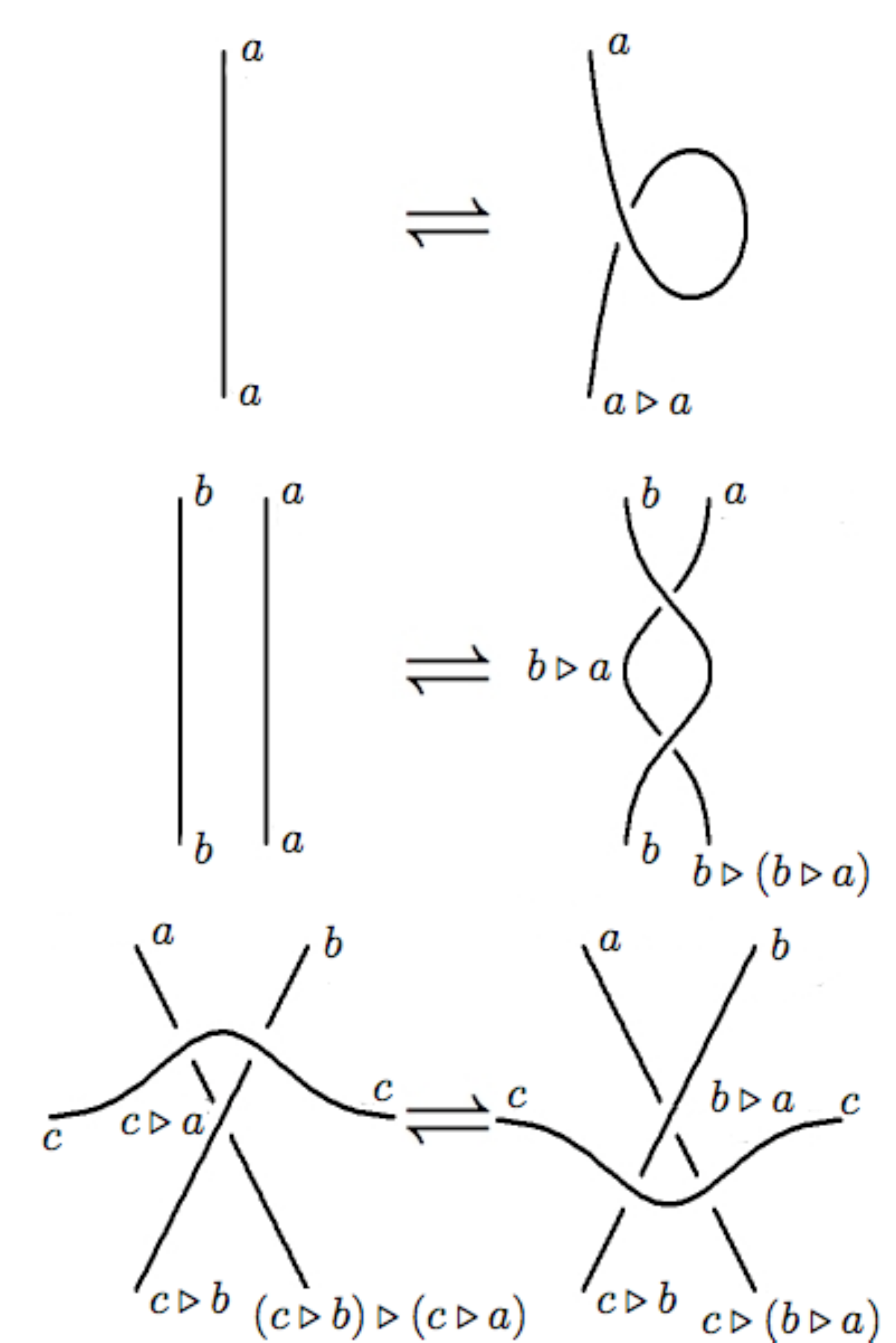}
\caption{The quandle axioms correspond to the Reidemeister moves}\label{fig:ReidemeisterColor}
\end{figure}

As it turns out, the number of colorings of a diagram for a given link by a given involutory quandle is independent of which diagram of the link we use.  Indeed, given a coloring of a link diagram, we get a unique coloring of any link diagram related to it by a Reidemeister move.  In fact, the three quandle axioms exactly correspond to the three Reidemeister moves, as indicated in figure \ref{fig:ReidemeisterColor}.

Thus we have the
\begin{Theorem}
For any involutory quandle $X$, the number of colorings of an unoriented link diagram by $X$ is an invariant of unoriented links.
\label{thm:basic}
\end{Theorem}

\subsection{The Fundamental Involutory Quandle}
Given an unoriented link diagram, we can define its fundamental involutory quandle\cite{WinkerThesis}.  This is a quandle which contains exactly the relations forced by the crossings in the diagram.  It is, in a sense, ``universal'' for colorings.

We generate a free quandle\cite{MR638121} on the set of arcs in the diagram $K$.  We then impose a relation for each crossing.  If generators $a$ and $c$ meet at the overcrossing generator $b$, we add the relation $c=b\triangleright a$.  Once these relations are added, the result is the fundamental involutory quandle $Q(K)$.

A coloring of the diagram $K$ by the quandle $X$ assigns to each arc of $K$ an element of $X$.  But these arcs are the generators of $Q(K)$.  Further, the relations defining $Q(K)$ are enforced by the definition of an $X$-coloring.  Thus an $X$-coloring of the link diagram $K$ is exactly the same as a quandle homomorphism $\hom_\Quan(Q(K),X)$.

When we apply a Reidemeister move to turn the diagram $K_1$ into the diagram $K_2$, the fundamental involutory quandle \textit{doesn't} stay the same.  The set of arcs in $K_2$ is not the same as the set of arcs in $K_1$, and there are different relations imposed by the different crossings.  However, we do have the

\begin{Theorem}
If link diagrams $K_1$ and $K_2$ are related by a Reidemeister move, then there is an isomorphism $Q(K_1)\cong Q(K_2)$.
\label{thm:FIQFunctorial}
\end{Theorem}
\begin{proof}
If we refer to figure \ref{fig:ReidemeisterColor} we can see the proof.  For example, let's say that $K_1$ is on the left side of a Reidemeister II move, while $K_2$ is on the right.

The labels in the middle row of figure \ref{fig:ReidemeisterColor} describe a coloring of $K_2$ using the quandle $Q(K_1)$, or equivalently a coloring of $K_1$ using the quandle $Q(K_2)$.  Thus wecan define two homomorphisms of quandles: $f\in\hom_\Quan(Q(K_2),Q(K_1))$ and $g\in\hom_\Quan(Q(K_1),Q(K_2))$.  These are clearly inverses of each other, establishing the isomorphism.
\end{proof}

In particular, this isomorphism gives a bijection between the sets of colorings $\hom_\Quan(Q(K_1),X)$ and $\hom_\Quan(Q(K_2),X)$, which reestablishes the invariance of coloring numbers.

It is important to note at this point that these sets of colorings are \emph{not} the same set.  They are merely \emph{isomorphic} as sets, rather than identical.  Therefore the set of colorings is not an invariant of the knot type.  Only its cardinality is invariant.  We must now lay out a language in which to talk about exactly these details.

\section{Categorification}
Categorification is, simply put
\begin{quote}
... the process of finding category-theoretic analogues of set-theoretic concepts by replacing sets with categories, functions with functors, and equations between functions by natural isomorphisms between functors, which in turn should satisfy certain equations of their own, called `coherence laws'.\cite{MR1664990}
\end{quote}
More to the point, we want to take things we'd called ``identical'' and see them as merely ``equivalent''.

In the case at hand, we're considering a knot to be an equivalence class of knot diagrams under the Reidemeister moves.  Instead, we'd like to think of link diagrams as the objects of a category $\KDiag$.  The morphisms will be sequences of Reidemeister moves.  Since any such move can be reversed, this category of link diagrams forms a groupoid.  Now we can recast theorem \ref{thm:basic} as follows:

\begin{Theorem}
For any involutory quandle $X$ we have a functor $\mathrm{Col}_X$ from the groupoid $\KDiag$ to the set of natural numbers, considered as a category with no non-identity morphisms.
\label{thm:categorified1}
\end{Theorem}
\begin{proof}
To any diagram we associate the number of $X$-colorings.  This defines the functor on objects.

Since every morphism is a composite of Reidemeister moves, we just need to define the functor on the Reidemeister moves to define it on all morphisms.  But we know that under a Reidemeister move the number of $X$-colorings remains the same, so to any move between two diagrams we can associate the identity morphism on the (common) number of colorings.
\end{proof}

We can also categorify the value of our invariant.  Instead of considering how many colorings a given diagram has, we should instead consider the set of colorings itself.  We further refine theorem \ref{thm:categorified1} to state:

\begin{Theorem}
For any involutory quandle $X$ we have a functor
\begin{equation*}
\mathrm{Col}_X:\KDiag\to\Set
\end{equation*}
which associates to any link diagram $K$ the set of $X$-colorings of $K$.
\label{thm:categorified2}
\end{Theorem}
\begin{proof}
Indeed, we can now see \ref{thm:FIQFunctorial} as asserting the functoriality of the fundamental involutory quandle construction.  That is, to a sequence of Reidemeister moves connecting two link diagrams we get an isomorphism of fundamental involutory quandles.  Then we can define $\mathrm{Col}_X(K)=\hom_\Quan(Q(K),X)$
\end{proof}

Thus a sequence of Reidemeister moves connecting two link diagrams gives an explicit bijection between the sets of $X$-colorings.  Since the sets are changing as we change the diagram, it no longer seems appropriate to call our functor a ``link invariant''.  Instead, we will make the following definition

\begin{Definition}
A \textit{link covariant} is a functor from the groupoid $\KDiag$ to any other category.  If the image of each morphism is an identity morphism, we call the functor a \textit{link invariant}.
\label{def:LinkCovariant}
\end{Definition}

Thus the fundamental involutory quandle of a knot diagram is a covariant, as is the set of $X$-colorings for any involutory quandle $X$.  Many other well-known ``invariants'' are actually covariants under this definition, like the knot group given by the Wirtinger presentation\cite{MR1211642}.

\section{Tangles}
\subsection{The 2-category of Tangles}
Now that we've categorified our link invariant, we have enough breathing room to truly extend its domain of definition.  Specifically, we want to color tangle diagrams.

Topologically, a tangle is like a knot or a link embedded in a cube, but we now allow arc components with their edges running to marked points on the top and bottom of the cube.  These tangles are known to form a monoidal category $\Tang$.  The objects of this category are the natural numbers, and a morphism from $m$ to $n$ is a tangle with $m$ points on the bottom of its cube, and $n$ endpoints on the top.

If we have a tangle from $n_1$ to $n_2$, and another tangle from $n_2$ to $n_3$, we can stack the second cube on top of the first and splice together the $n_2$ endpoints in the middle.  This defines our composition.  The monoidal product of two objects is their sum as natural numbers, while the monoidal product of two tangles is given by stacking their cubes side-by-side.

Just as for knots and links, tangles can be described by tangle diagrams.  Ambient isotopies of tangles are again equivalent to sequences of Reidemeister moves.  This leads to a well-known presentation of $\Tang$ as a monoidal category\cite{MR1020583}:

\begin{Theorem}
The category $\Tang$ of tangle diagrams is generated by the tangle diagrams $\{X^+,X^-,\cup,\cap\}$ with relations
\begin{list}{}{}
\item{$T_0$.} $(\cup\otimes I_1)\circ(I_1\otimes\cap) = I_1 = (I_1\otimes\cup)\circ(\cap\otimes I_1)$
\item{$T_0'$.} $(I_1\otimes\cup)\circ(X^\pm\otimes I_1) = (\cup\otimes I_1)\circ(I_1\otimes X^\mp)$
\item{$T_1$.} $\cup\circ X^\pm = \cup$
\item{$T_2$.} $X^\pm\circ X^\mp = I_2$
\item{$T_3$.} $(X^+\otimes I_1)\circ(I_1\otimes X^+)\circ(X^+\otimes I_1)=(I_1\otimes X^+)\circ(X^+\otimes I_1)\circ(I_1\otimes X^+)$
\end{list}
\end{Theorem}

We read the generator $X^+$ as a right-handed crossing, $X^-$ as a left-handed crossing, $\cup$ as a local minimum in the tangle diagram, and $\cap$ as a local maximum.  The relations $T_1$, $T_2$, and $T_3$ then encode the three Reidemeister moves, while $T_0$ and $T_0'$ handle the interaction of local maxima and minima with each other and with crossings.

As we did before, let's categorify this picture.  Instead of identifying two tangle diagrams if they are related by a Reidemeister move (or one of the new ``topological'' tangle moves), let's jut consider them to be equivalent.

That is, we consider a (strict) monoidal 2-category whose objects are again the natural numbers, and whose morphisms are built from compositions and monoidal products of the four generating tangles.  Now instead of imposing the five relations, we add 2-isomorphisms to relate any tangle diagrams that would be identified by the relations.  It is this 2-category that we will refer to as $\Tang$.

In analogy with definition \ref{def:LinkCovariant} for links, we introduce the following
\begin{Definition}
A \textit{tangle covariant} is a monoidal 2-functor from the monoidal 2-category $\Tang$ to any other 2-category.  If the image of each 2-morphism is an identity 2-morphism, we call the functor a \textit{tangle invariant}.
\label{def:TangleCovariant}
\end{Definition}

The straightforward approach now is to define a coloring of an unoriented tangle diagram by an involutory quandle $X$ exactly as we did for link diagrams.  We assign an element of $X$ to each arc and subject these assignments to restrictions at crossings just as before.  This indeed gives a set of $X$-colorings, but there is no way to compose two of these sets as morphisms in some category.  We need to extend our na\:ive notion of the set of tangle colorings and give it ``handles'' that we can use to compose them.

\section{Spans}
\subsection{The 2-category of spans}
Given a category $\mathcal{C}$ with pullbacks we define the 2-category $\Span(\mathcal{C})$ of spans on $\mathcal{C}$.  It will have the same objects as $\mathcal{C}$.  A morphism $f:A\to B$ in $\Span(\mathcal{C})$ will be a ``span'' in $\mathcal{C}$: an object $F$ and a pair of morphisms in $\mathcal{C}$: $A\xleftarrow{f_l}F\xrightarrow{f_r}B$.  Then, given spans $f=A\xleftarrow{f_l}F\xrightarrow{f_r}B$ and $g=A\xleftarrow{g_l}G\xrightarrow{g_r}B$, a 2-morphism $\phi:f\Rightarrow g$ is an arrow $\phi:F\to G$ so that the following diagram commutes:
\begin{diagram}
A&\lTo^{f_l}&F\\
\uTo^{g_l}&\ldTo~{\phi}&\dTo_{f_r}\\
G&\rTo_{g_r}&B\\
\end{diagram}

The ``vertical'' composition of 2-morphisms is straightforward.  The composition of morphisms (and the ``horizontal'' composition of 2-morphisms) invokes the pullbacks we assumed $\mathcal{C}$ to have.  If we have spans $f=A\xleftarrow{f_l}F\xrightarrow{f_r}B$ and $g=B\xleftarrow{g_l}G\xrightarrow{g_r}C$ we form their composite by pulling back the square in the diagram
\begin{diagram}
F\circ G\SEpbk&\rTo&G&\rTo^{g_r}&C\\
\dTo&&\dTo^{g_l}&&\\
F&\rTo^{f_r}&B&&\\
\dTo^{f_l}&&&&\\
A&&&&\\
\end{diagram}
This composition is not quite associative, but it's easily verified to be associative up to a unique 2-isomorphism, which gives the associator for the 2-category.

There are a few facts about the span construction which will be useful to us.\cite{ArmstrongThesis}
\begin{Theorem}
Given categories $\mathcal{C}$ and $\mathcal{D}$ with pullbacks and a functor $F:\mathcal{C}\to\mathcal{D}$ preserving them, there is a 2-functor $\Span(F):\Span(\mathcal{C})\to\Span(\mathcal{D})$ defined by applying $F$ to all parts of a span diagram.
\end{Theorem}
\begin{Theorem}
If $\mathcal{C}$ is a monoidal category such that the monoidal product preserves pullbacks, then $\Span(\mathcal{C})$ is a monoidal 2-category.
\end{Theorem}

Dually, given a category $\mathcal{C}$ with pushouts we can define the 2-category $\CoSpan(\mathcal{C})$ of cospans.  A cospan diagram is like a span diagram, but with the arrows pointing in instead of out, and we compose them by pushing out a square rather than pulling back, but otherwise everything we've said about spans holds for cospans.

\subsection{Coloring Spans}
The category $\Set$ of sets has fibered products, which act as pullbacks, and so we have a 2-category $\Span(\Set)$.  The two functions out to the side of the central set in a span will provide us with exactly the handles we need to compose sets of colorings.

Now we can extend theorem \ref{thm:categorified2} to:
\begin{Theorem}
For any involutory quandle $X$ we have a 2-functor
\begin{equation*}
\mathrm{Col}_X:\Tang\to\Span(\Set)
\end{equation*}

On an object $n$ of $\Tang$ we define $\mathrm{Col}_X(n)=X^n$ the set of $n$-tuples of elements of $X$.

For a tangle diagram $T:m\to n$ from $m$ free ends to $n$ free ends we define the span
\begin{equation*}
X^m\leftarrow\mathrm{Col}_X(T)\rightarrow X^n
\end{equation*}
where the arrow on the left is the function sending a coloring of $T$ to the coloring it induces on the lower endpoints of the tangle, and the one on the right is the similar function for the upper endpoints.

The 2-functor is defined on 2-morphisms by the diagrams in figure \ref{fig:ReidemeisterColor}, as in theorem \ref{thm:categorified2}.
\label{thm:spanofsets}
\end{Theorem}
\begin{proof}
The main thing to check here is that composition of coloring spans really does reflect composition of tangles.  But given a composite tangle $T_1\circ T_2$, a coloring in $\mathrm{Col}_X(T_1\circ T_2)$ is exactly a coloring of $T_1$ and a coloring of $T_2$ that agree on the endpoints we splice together to compose the tangles.  This is exactly the definition of the fibered product
\begin{diagram}
\mathrm{Col}_X(T_1\circ T_2)&\rTo&\mathrm{Col}_X(T_2)&\rTo&X^p\\
\dTo&&\dTo&&\\
\mathrm{Col}_X(T_1)&\rTo&X^n&&\\
\dTo&&&&\\
X^m&&&&\\
\end{diagram}
\end{proof}

Notice what happens to this picture when we consider a link as a tangle from $0$ to $0$.  Both sides of the span become empty products -- singletons -- and the functions in the span become trivial.  What remains is the old set of link colorings.

\subsection{The Fundamenal Involutory Quandle Cospan}
Earlier we identified the fundamental involutory quandle $Q(K)$ of a link diagram $K$ as the quandle that captures coloring numbers for all involutory quandles $X$:
\begin{equation*}
\mathrm{Col}_X(K)\cong\hom_\Quan(Q(K),X)
\end{equation*}

The same construction can give us a quandle $Q(T)$ from a tangle diagram $T$, which then gives us the set of $X$-colorings of $T$.  Can we get the sides of our span as well?

Indeed, the free quandle on $n$ generators $Q_n$ satisfies $\hom_\Quan(Q_n,X)=X^n$.  We can choose these generators to be a collection of free ends of our tangle diagram, and the inclusion of those ends into the whole diagram gives us a homomorphism $Q_n\to Q(T)$.

\begin{Theorem}
There is a 2-functor extending the fundamental involutory quandle to tangles:
\begin{equation*}
Q:\Tang\to\CoSpan(\Quan)
\end{equation*}

On an object $n$ in $\Tang$ we define $Q(n)=Q_n$, the free quandle on $n$ generators.

For a tangle diagram $T:m\to n$ from $m$ free ends to $n$ free ends we let $Q_m$ be the free quandle on the incoming ends and $Q_n$ be the free quandle on the outgoing ends.  We define the cospan
\begin{equation*}
Q_m\rightarrow Q(T)\leftarrow Q_n
\end{equation*}
where the arrows are the quandle homomorphisms induced by including the endpoints into the tangle diagrams.

For a 2-morphism $\phi$ we define $Q(\phi)$ by referring to figure \ref{fig:ReidemeisterColor}, as in theorem \ref{thm:FIQFunctorial}.
\end{Theorem}
\begin{proof}
Again, the meat of the proof is in showing that composition of tangles really does correspond to a pushout in $\Quan$.

Composition of tangle diagrams $T_1$ and $T_2$ consists of laying down both diagrams and joining some arcs from $T_1$ to arcs from $T_2$, as determined by the lineup of the endpoints.  But matching endpoints corresponds to adding relations saying that the image of a generator of $Q_n$ in $Q(T_1)$ equals its image as a generator in $Q(T_2)$.  This amalgamated free product is exactly the pushout construction in $\Quan$.
\end{proof}

Again, if we consider a link as a tangle from $0$ to $0$, the free quandle on zero generators is trivial, as are all homomorphisms from it.  The only nontrivial information in this cospan is the old fundamental involutory quandle of the link.

Now we can use this fundamental involutory quandle cospan to recover the coloring spans.  The contravariant hom-functor $\hom_\Quan(\underline{\hphantom{X}},X)$ automatically takes all colimits to limits, so in particular it preserves pullbacks as a functor $\Quan^{\mathrm{op}}\to\Set$.

\begin{Theorem}
The coloring span 2-functor $\mathrm{Col}_X$ factors as the composition of the span of the hom-functor $\Span(\hom_\Quan(\underline{\hphantom{X}},X))$ and the fundamental involutory quandle 2-functor $Q$.
\label{thm:coloringFactors}
\end{Theorem}

\subsection{Monoidal structure}
All of the 2-categories considered above also carry monoidal structures, and all the 2-functors preserve them.  This allows us to obtain tangle covariants, and to decategorify them to tangle invariants.

The category $\Set$ has all finite products, so it has the Cartesian monoidal structure.  The direct product of sets preserves pullbacks, so $\Span(\Set)$ is a monoidal 2-category.

Similarly, $\Quan$ has finite coproducts given by the free product of quandles, or equivalently by the pushout over the free quandle on zero generators.  These coproducts preserve pushouts, so $\CoSpan(\Quan)$ is a monoidal 2-category.

\begin{Theorem}
The induced 2-functor
\begin{equation*}
\Span(\hom_\Quan(\underline{\hphantom{X}},X)):\CoSpan(\Quan)\to\Span(\Set)
\end{equation*}
is monoidal.
\end{Theorem}
\begin{proof}
This is a straightforward consequence of the fact that the hom-functor $\hom_\Quan(\underline{\hphantom{X}},X):\Quan^{\mathrm{op}}\to\Set$ preserves products.
\end{proof}

\begin{Theorem}
The fundamental involutory quandle cospan 2-functor
\begin{equation*}
Q:\Tang\to\CoSpan(\Quan)
\end{equation*}
is monoidal.
\end{Theorem}
\begin{proof}
Given two tangles $T_1:m_1\to n_1$ and $T_2:m_2\to n_2$ we form their monoidal product $T_1\otimes T_2$ by laying them side-by-side.  When we calculate the fundamental involutory quandle of this diagram, we just use all the generators and relations that come from each of $T_1$ and $T_2$, and none of them interact with each other.  Thus the quandle of $T_1\otimes T_2$ is the free product of the quandles of $T_1$ and $T_2$.  Similarly at the ends, $Q_{m_1+m_2}$ is the free product of $Q_{m_1}$ and $Q_{m_2}$, and $Q_{n_1+n_2}$ is the free product of $Q_{n_1}$ and $Q_{n_2}$.  So the monoidal product of tangles corresponds under $Q$ to taking free products of cospan diagrams.  But this is just the induced monoidal structure on $\CoSpan(\Quan)$.
\end{proof}

\begin{Theorem}
For any involutory quandle $X$ the coloring span 2-functor
\begin{equation*}
\mathrm{Col}_X:\Tang\to\Span(\Set)
\end{equation*}
is monoidal.
\end{Theorem}
\begin{proof}
This is an immediate corollary of the preceding theorems and theorem \ref{thm:coloringFactors}
\end{proof}

\section{Decategorifying}

\subsection{Coloring Matrices}
When we decategorify a coloring set we get a coloring number.  What happens when we decategorify a coloring span?

A 2-isomorphism in the 2-category $\Span(\Set)$ is a bijection $\phi:F\rightarrow G$ in diagram
\begin{diagram}
A&\lTo^{f_l}&F\\
\uTo^{g_l}&\ldTo~{\phi}&\dTo_{f_r}\\
G&\rTo_{g_r}&B\\
\end{diagram}
The span functions $f_l$ and $f_r$ partition $F$ into its ``double preimages''
\begin{equation*}
F=\bigcup_{\substack{a\in A\\b\in B}}F_{a,b}\qquad\qquad F_{a,b}=\{x\in F|f_l(x)=a, f_r(x)=b\}
\end{equation*}

Similarly, the functions $g_l$ and $g_r$ partition $G$ into its double preimages $G_{a,b}$.  Then for the diagram above to commute the function $\phi$ must decompose into functions $\phi_{a,b}:F_{a,b}\to G_{a,b}$.  And then for $\phi$ to be a bijection, each of the $\phi_{a,b}$ must be a bijection.

So when we identify isomorphic spans of sets, we retain only the cardinality of each of the double preimages.  We are left with a matrix of cardinal numbers indexed by the set $A$ on the one side and the set $B$ on the other.

For a coloring span, these index sets are the colorings of the endpoints.  Thus when we decategorify a coloring span we get a matrix $\mathrm{Col}_X(T)$ indexed by colorings of the endpoints of the tangle.  The entry $\mathrm{Col}_X(T)_{\mu\nu}$ is the number of colorings of the diagram $T$ that agree with the coloring $\mu$ on the incoming ends and with the coloring $\nu$ on the outgoing ends.

This interpretation as matrices is compatible with matrix multiplication.  That is, given tangle diagrams $T_1:m\to l$ and $T_2:l\to n$, the number of colorings $\mathrm{Col}_X(T_1\circ T_2)_{\mu\nu}$ agreeing with the colorings $\mu$ and $\nu$ on the ends can be calculated as a sum of products of coloring numbers:
\begin{equation*}
\mathrm{Col}_X(T_1\circ T_2)_{\mu\nu}=\sum_{\lambda\in X^l}\mathrm{Col}_X(T_1)_{\mu\lambda}\mathrm{Col}_X(T_2)_{\lambda\nu}
\end{equation*}

Decategorification also plays nice with the monoidal structure on spans induced by the product of sets.  Take two diagrams $T_1:m_1\to n_1$ and $T_2:m_2\to n_2$.  A coloring $\mu_1$ of the incoming ends of $T_1$ and a coloring $\mu_2$ of the incoming ends of $T_2$ combine to give a coloring $(\mu_1,\mu_2)\in X^{m_1+m_2}$ of the incoming ends of $T_1\otimes T_2$.  Similarly, we can combine colorings of the outgoing strands of each diagram to get a coloring $(\nu_1,\nu_2)\in X^{n_1+n_2}$ of the outgoing strands of $T_1\otimes T_2$.  Every coloring of the incoming or outgoing strands arises in this manner.

Now when we count the colorings of $T_1\otimes T_2$ compatible with a given coloring of the incoming and outgoing ends, we find
\begin{align*}
\mathrm{Col}_X(T_1\otimes T_2)_{(\mu_1,\mu_2)(\nu_1,\nu_2)}&=\mathrm{Col}_X(T_1)_{\mu_1\nu_1}\mathrm{Col}_X(T_2)_{\mu_2\nu_2}\\
&=\left(\mathrm{Col}_X(T_1)\boxtimes\mathrm{Col}_X(T_2)\right)_{(\mu_1,\mu_2)(\nu_1,\nu_2)}\\
\end{align*}
This follows since a coloring of $T_1\otimes T_2$ is simply a coloring of each of $T_1$ and $T_2$ with no particular relation between them.  This shows that the coloring matrix for the monoidal product $T_1\otimes T_2$ is the Kronecker product of the coloring matrices for $T_1$ and $T_2$.

\begin{Theorem}
For any finite involutory quandle $X$, there is a monoidal 2-functor
\begin{equation*}
\mathrm{Col}_X:\Tang\to\Mat(\mathbb{N})
\end{equation*}
where the target category is that of matrices with natural number entries, and with identity 2-morphisms added.
\end{Theorem}
\begin{proof}
If we pick $d$ to be the cardinality of $X$, then there are exactly $d^n$ colorings of a collection of $n$ endpoints in a tangle.  We thus set $\mathrm{Col}_X(n)=d^n$ on objects.

We already have a coloring span of sets for every tangle.  Even if we disregard the coloring relations at crossings, we can only pick one color from $X$ for each arc in the diagram, and so the sets in the coloring span are finite.  Taking cardinalities, we get a matrix of natural numbers.  As described above, this assignment of a coloring matrix to a tangle preserves the composition and monoidal structure.

Finally, if we have a 2-morphism $\phi:T_1\Rightarrow T_2$ in $\Tang$ we know that the coloring matrices for $T_1$ and $T_2$ will be the same, so we can pick $\mathrm{Col}_X(\phi)$ to be the identity 2-morphism on that matrix.
\end{proof}

Since every 2-morphism becomes an identity 2-morphism under this functor, we have a tangle invariant.

In particular, when we consider a link $L$ as a tangle from $0$ to $0$, we can find the $1\times1$ matrix $\mathrm{Col}_X(L)$.  The single entry in this matrix is the number of $X$-colorings of the link $L$.

Instead of restricting our attention to links, we may instead consider any $n$-strand braid as a tangle from $n$ to $n$.  In this case we find a matrix representation $\mathrm{Col}_X$ of each braid group $B_n$.

\subsection{Computation}
It turns out that not only do we have a tangle invariant in our coloring matrices, we have a straightforward way of computing them.  The category of tangles was given by generators and relations.  Thus we can calculate the coloring matrix of each generating tangle by hand, and then assemble the coloring matrix using matrix multiplications and Kronecker products.

The matrix for each generating tangle is straightforward to work out.  The right-handed crossing, for instance, takes a pair of colors for each index.  The entry $\mathrm{Col}_X(X^+)_{(a,b)(c,d)}$ will be $1$ if $a=d$ and $c=a\triangleright b$, and $0$ otherwise.  As an example, figure \ref{fig:DZ3colors} shows all the coloring matrices of the generating tangles for the quandle $D(\mathbb{Z}_3)$.

\begin{figure}[hbtp]
  \begin{align*}
    \mathrm{Col}_{D(\mathbb{Z}_3)}(\cup)&=\begin{pmatrix}1\\0\\0\\0\\1\\0\\0\\0\\1\end{pmatrix}\\
    \mathrm{Col}_{D(\mathbb{Z}_3)}(\cap)&=\begin{pmatrix}1&0&0&0&1&0&0&0&1\end{pmatrix}\\
    \mathrm{Col}_{D(\mathbb{Z}_3)}(X^+)&=\begin{pmatrix}1&0&0&0&0&0&0&0&0\\0&0&0&0&0&0&1&0&0\\0&0&0&1&0&0&0&0&0\\0&0&0&0&0&0&0&1&0\\0&0&0&0&1&0&0&0&0\\0&1&0&0&0&0&0&0&0\\0&0&0&0&0&1&0&0&0\\0&0&1&0&0&0&0&0&0\\0&0&0&0&0&0&0&0&1\end{pmatrix}\\
    \mathrm{Col}_{D(\mathbb{Z}_3)}(X^-)&=\begin{pmatrix}1&0&0&0&0&0&0&0&0\\0&0&0&0&0&1&0&0&0\\0&0&0&0&0&0&0&1&0\\0&0&1&0&0&0&0&0&0\\0&0&0&0&1&0&0&0&0\\0&0&0&0&0&0&1&0&0\\0&1&0&0&0&0&0&0&0\\0&0&0&1&0&0&0&0&0\\0&0&0&0&0&0&0&0&1\end{pmatrix}\\  \end{align*}
\caption{$D(\mathbb{Z}_3)$-coloring matrices for the generators of $\Tang$}
\label{fig:DZ3colors}
\end{figure}

Computations with these matrices may be tedious by hand, but they are easily programmed into a computer.

\bibliographystyle{amsplain}
\bibliography{../biblio}
\end{document}